\documentclass[a4paper,10pt]{amsart}
\usepackage{amsmath,amssymb,amsfonts,amsthm}
\usepackage{latexsym,mathrsfs}
\usepackage[all]{xy}
\usepackage{longtable}
\usepackage{supertabular}
\input xypic
\theoremstyle{plain}
\newtheorem{theorem}[subsection]{{\bf Theorem}}
\newtheorem*{theorem*}{{\bf Theorem}}

\newtheorem*{corollary*}{{\bf Corollary}}

\newtheorem{question}[subsection]{{\bf Question}}

\theoremstyle{definition}

\theoremstyle{remark}

\numberwithin{equation}{subsection}
\DeclareMathOperator{\HH}{H}

\DeclareMathOperator{\B}{B}
\DeclareMathOperator{\Bt}{\tilde{B}}

\DeclareMathOperator{\Hom}{Hom}

\DeclareMathOperator{\res}{res}

\newcommand{\QZ}{\mathbb{Q}/\mathbb{Z}}

\newcommand{\cwedge}{\curlywedge}

\begin{document}
\title[Unramified Brauer groups]
{Unramified Brauer groups and isoclinism}
\author{Primo\v z Moravec}
\address{{
Department of Mathematics \\
University of Ljubljana \\
Jadranska 21 \\
1000 Ljubljana \\
Slovenia}}
\email{primoz.moravec@fmf.uni-lj.si}
\thanks{}
\date{\today}
%%%%%%%%%%%%%%%%%%%%%%%%%%%%%%%%%%%%%%%%%%%%%%%%%%%%%%%%%%%%%%%%%%%%%
\begin{abstract}
\noindent
We show that if $G_1$ and $G_2$ are isoclinic groups, then their Bogomolov
multipliers are isomorphic.
\end{abstract}
%%%%%%%%%%%%%%%%%%%%%%%%%%%%%%%%%%%%%%%%%%%%%%%%%%%%%%%%%%%%%%%%%%%%%%%
\maketitle
%%%%%%%%%%%%%%%%%%%%%%%%%%%%%%%%%%%%%%%%%%%%%%%%%%%%%%%%%%%%%%%%%%%%%%%
\section{Introduction}
\label{s:intro}

\noindent
Let $G$ be a finite group and $V$ a faithful representation of $G$ over $\mathbb{C}$. Then there is a natural action of $G$ upon
the field of rational functions $\mathbb{C}(V)$. Noether's problem \cite{Noe16} asks as to whether 
the field of $G$-invariant functions $\mathbb{C}(V)^G$ is purely transcendental over $\mathbb{C}$, i.e.,
whether the quotient space $V/G$ is {\it rational}. A question related to the above mentioned
is whether $V/G$ is {\it stably rational}, that is,
whether there exist independent variables $x_1,\ldots ,x_r$ such that $\mathbb{C}(V)^G(x_1,\ldots ,x_r)$ becomes a pure
transcendental extension of $\mathbb{C}$. This problem has close connection
with L\"{u}roth's problem \cite{Saf91} and the inverse Galois problem \cite{Swa83,Sal84}.
By Hilbert's Theorem 90 stable rationality of $V/G$
does not depend upon the choice of $V$, but only on the group $G$. Saltman \cite{Sal84} found examples of groups $G$ of order $p^9$ such
that $V/G$ is not stably rational over $\mathbb{C}$. His main method was application of the unramified cohomology group
$\HH^2_{\rm nr}(\mathbb{C}(V)^G,\QZ )$ as an obstruction. 
A version of this invariant had been used before by Artin and Mumford \cite{Art72} who constructed unirational varieties over
$\mathbb{C}$ that were not rational.
Bogomolov \cite{Bog88} further explored this cohomology group. He proved 
that $\HH^2_{\rm nr}(\mathbb{C}(V)^G,\QZ )$ is canonically isomorphic to 
\begin{equation}
\label{eq:b0def}
\B_0(G)=\bigcap _{\tiny\begin{matrix}A\le G,\\ A \hbox{ abelian}\end{matrix}} \ker \res ^G_A,
\end{equation}
where $\res ^G_A:\HH^2(G,\QZ )\to \HH^2(A,\QZ )$ is the usual cohomological restriction map.
The group $\B_0(G)$ is a subgroup of the {\it Schur multiplier} $\HH ^2(G,\QZ )$ of $G$. Kunyavski\u\i ~\cite{Kun08} coined the term the
{\it Bogomolov multiplier} of $G$ for the group $\B_0(G)$.

We recently proved \cite{Mor11} that $\B_0(G)$ is naturally isomorphic to $\Hom (\Bt_0(G),\QZ)$, where $\Bt_0(G)$ is the kernel of the commutator
map $G\cwedge G\to [G,G]$, and $G\cwedge G$ is a quotient of the {\it non-abelian exterior square} of $G$ (see Section \ref{s:proof} for further
details). This description of $\B_0(G)$ is purely combinatorial, and allows for efficient computations of $\B_0(G)$, and a Hopf formula for
$\B_0(G)$. We also note here that the group $\Bt_0(G)$ can be defined for any (possibly infinite) group $G$. 

Recently, Hoshi, Kang, and Kunyavski\u\i ~\cite{Hos12} classified all groups of order $p^5$ with non-trivial Bogomolov multiplier; another
classification was found in \cite{Mor12}. It turns out that only examples of such groups appear within the same isoclinism family, where
isoclinism is the notion defined by P. Hall \cite{Hal40}.
The following question was posed in \cite{Hos12}:

\begin{question}[\cite{Hos12}]
\label{q:iso}
Let $G_1$ and $G_2$ be isoclinic $p$-groups. Is it true that the fields $k(G_1)$ and $k(G_2)$ are stably isomorphic, or at least,
that $\B_0(G_1)$ is isomorphic to $\B_0(G_2)$?
\end{question}

The purpose of this note is to answer the second part of the above question in the affirmative:

\begin{theorem}
\label{t:iso}
Let $G_1$ and $G_2$ be isoclinic groups. Then $\Bt_0(G_1)\cong \Bt_0(G_2)$. In particular, if $G_1$ and $G_2$ are finite,
then $\B_0(G_1)$ is isomorphic to $\B_0(G_2)$.
\end{theorem}

%%%%%%%%%%%%%%%%%%%%%%%%%%%%%%%%%%%%%%%%%%%%%%%%%%%%%%%%%%%%%%%%%%%%%%%%%%%%%%%%%%%%%%%%%%%%%%%%%%%%%%%%%%%%%%%%%%
\section{Proof of Theorem \ref{t:iso}}
\label{s:proof}

\noindent
We first recall the definition of $G\cwedge G$ from \cite{Mor11}.
For $x,y\in G$ we write ${}^xy=xyx^{-1}$ and $[x,y]=xyx^{-1}y^{-1}$.
Let $G$ be any group (possibly infinite). We form the group $G\cwedge G$, generated by the symbols
$m\cwedge n$, where $m,n\in G$, subject to the following relations:
\begin{align}
mm'\cwedge n &= ({}^mm'\cwedge {}^mn)(m\cwedge n),\notag \\
m\cwedge nn' &= (m\cwedge n)({}^nm\cwedge {}^nn'), \label{eq:cwed}\\
x\cwedge y & = 1,\notag
\end{align}
for all $m,m',n,n'\in G$, and all $x,y\in G$ with $[x,y]=1$. The group $G\cwedge G$ is a quotient of the non-abelian exterior square
$G\wedge G$ of $G$ defined by Miller \cite{Mil52}. There is a surjective homomorphism $\kappa :G\cwedge G\to [G,G]$ defined by
$\kappa (x\cwedge y)=[x,y]$ for all $x,y\in G$. Denote $\Bt _0(G)=\ker\kappa$. By \cite{Mor11} we have the following:

\begin{theorem}[\cite{Mor11}]
\label{t:mor}
Let $G$ be a finite group. Then $\B_0(G)$ is naturally isomorphic to $\Hom(\Bt_0(G),\QZ)$, and thus
$\B_0(G)\cong \Bt_0(G)$.
\end{theorem} 

Let $L$ be a group. A function $\phi : G \times G \to L$ is called a {\it $\Bt_0$-pairing} if for all $m, m',n, n' \in G$,
and for all $x,y\in G$ with $[x,y]=1$,
\begin{align*}
\phi (mm', n) &= \phi ({}^mm', {}^mn)\phi (m,n),\\
\phi (m,nn') &= \phi (m,n)\phi ({}^nm, {}^nn'),\\
\phi (x,y) &=1.
\end{align*}
Clearly a $\Bt_0$-pairing $\phi$ determines a unique homomorphism of groups 
$\phi ^*: G\cwedge G\to L$ such that $\phi ^*(m\cwedge n)=\phi (m,n)$
for all $m,n\in G$.

We now turn to the proof of Theorem \ref{t:iso}. Let $G_1$ and $G_2$ are isoclinic groups, and denote 
$Z_1=Z(G_1)$, $Z_2=Z(G_2)$. By definition \cite{Hal40}, there exist
isomorphisms $\alpha :G_1/Z_1\to G_2/Z_2$ and $\beta :[G_1,G_1]\to [G_2,G_2]$ such that if $\alpha(a_1Z_1)=a_2Z_2$ and
$\alpha (b_1Z_1)=b_2Z_2$, then $\beta ([a_1,b_1])=[a_2,b_2]$ for all $a_1,b_1\in G_1$. Define a map $\phi :G_1\times G_1\to G_2\cwedge G_2$
by $\phi (a_1,b_1)=a_2\cwedge b_2$, where $a_i,b_i$ are as above. To see that this is well defined, suppose that
$\alpha (a_1Z_1)=a_2Z_2=\bar{a_2}Z_2$ and $\alpha (b_1Z_1)=b_2Z_2=\bar{b_2}Z_2$. Then we can write $\bar{a_2}=a_2z$ and $\bar{b_2}=b_2w$ for some $w,z\in Z_2$. By definition
of $G_2\cwedge G_2$ we have that $\bar{a_2}\cwedge \bar{b_2}=a_2\cwedge b_2$, hence $\phi$ is well defined. 

Suppose that $a_1,b_1\in G_1$ commute, and let $a_2,b_2\in G_2$ be as above. By definition, $[a_2,b_2]=\beta ([a_1,b_1])=1$, hence $a_2\cwedge b_2=1$.
This, and the relations of $G_2\cwedge G_2$, ensure that $\phi$ is a $\Bt_0$-pairing. Thus $\phi$ induces a homomorphism $\gamma :G_1\cwedge G_1\to
G_2\cwedge G_2$ such that $\gamma (a_1\cwedge b_1)=a_2\cwedge b_2$ for all $a_1,b_1\in G_1$. By symmetry there exists a homomorphism $\delta :G_2\cwedge G_2\to G_1\cwedge G_1$ defined via $\alpha ^{-1}$. It is straightforward to see that $\delta$ is the inverse of $\gamma$, hence $\gamma$ is an isomorphism.

Let $\kappa _1:G_1\cwedge G_1\to [G_1,G_1]$ and $\kappa _2:G_2\cwedge G_2\to [G_2,G_2]$ be the commutator maps. Since $\beta\kappa _1(a_1\cwedge b_1)=\beta ([a_1,b_1])=[a_2,b_2]=\kappa _2\gamma (a_1\cwedge b_1)$, we have the following commutative diagram with exact rows:

$$
\xymatrix{
0\ar[r] & \Bt_0(G_1)\ar[r]\ar[d]_{\tilde{\gamma}} & G_1\cwedge G_1 \ar[r]^{\kappa _1}\ar[d]_{\gamma} & [G_1,G_1]
 \ar[r]\ar[d]_{\beta} & 0\\
0\ar[r] & \Bt_0(G_2)\ar[r] & G_2\cwedge G_2\ar[r]^{\kappa _2} & [G_2,G_2]\ar[r] & 0
}
$$
Here $\tilde{\gamma}$ is the restriction of $\gamma$ to $\Bt_0(G_1)$. Since $\beta$ and $\gamma$ are isomorphisms, so is $\tilde{\gamma}$.
This concludes the proof.
%$$
%\xymatrix{
%0\ar[r] & \Bt _0(G_1)\ar[r]\ar[d]_{\gamma \mid _{\Bt_ 0(G_1)}} & G_1\cwedge G_1\ar[r]^{\kappa _1}\ar[d]_{\alpha} & [G_1,G_1]
% \ar[r]\ar[d]_{\beta} & 0\\
%0\ar[r] & \Bt _0(G_2)\ar[r] & G_2\cwedge G_2\ar[r]^{\kappa _2} & \[G_2,G_2]\ar[r] & 0
%}
%$$
%%%%%%%%%%%%%%%%%%%%%%%%%%%%%%%%%%%%%%%%%%%%%%%%%%%%%%%%%%%%%%%%%%%%%%%%%%%%%%%%%%%%%%%%%%%%%%%%%%%%%%%%%%%%%%%%%%%%%%%


\begin{thebibliography}{99}
%
\bibitem{Art72}
M. Artin, and D. Mumford,
{\it Some elementary examples of unirational varieties which are not rational},
Proc. London. Math. Soc. (3) {\bf 25} (1972), 75--95.
% 
%\bibitem{Bla79}
%N. Blackburn, and L. Evens,
%{\it Schur multipliers of $p$-groups},
%J. Reine Angew. Math. {\bf 309} (1979), 100--113. 
%
%\bibitem{Bly09}
%R. D. Blyth, and R. F. Morse, {\it Computing the nonabelian tensor square of polycyclic groups},
%J. Algebra {\bf 321} (2009), no. 8, 2139--2148.
%
\bibitem{Bog88}
F. A. Bogomolov, {\it The Brauer group of quotient spaces by linear group actions},
Izv. Akad. Nauk SSSR Ser. Mat {\bf 51} (1987), no. 3, 485--516.
%
%\bibitem{Bog04}
%F. A. Bogomolov, J. Maciel, T. Petrov,
%{\it Unramified Brauer groups of finite simple groups of type $A_\ell$},
%Amer. J. Math. {\bf 126} (2004), 935--949.
%
%\bibitem{Bro82}
%K. S. Brown, {\it Cohomology of groups}, Springer-Verlag, New York, Heidelberg, Berlin, 1982.
%
%\bibitem{Bro84}
%R. Brown, and J.-L. Loday, 
%{\it Excision homotopique en basse dimension},
%C. R. Acad. Sci. Paris S\'{e}r. Math. {\bf 289} (1984), 353--356.
%
%\bibitem{Bro87}
%R. Brown, and J.-L. Loday, 
%{\it Van Kampen theorems for diagrams of spaces},
%Topology {\bf 26} (1987), no. 3, 311--335.
%
%\bibitem{Chu09}
%H. Chu, S. Hu., M. Kang, and B. E. Kunyavski\u\i, {\it Noether's problem and the unramified Brauer group for groups of order 64},
%Int. Math. Res. Not. IMRN {\bf 12} (2010), 2329--2366.
%
%\bibitem{Chu08}
%H. Chu, S. Hu., M. Kang, and Y. G. Prokhorov, {\it Noether's problem for groups of order 32},
%J. Algebra {\bf 320} (2008), 3022--3035.
%
%\bibitem{Col89}
%J.-L. Colliot-Th\'{e}l\`{e}ne, and M. Ojanguren, {\it Vari\'{e}t\'{e}s unirationelles non rationelles: au-del\`{a} de l'exemple d'Artin et Mumford},
%Invent. Math. {\bf 97} (1989), 141--158.
%
%\bibitem{Den73}
%R. K. Dennis, and M. R. Stein, {\it The functor $\K_2$: a survey of computations and problems},
%Algebraic $K$-theory, II: ``Classical'' algebraic K-theory and connections with arithmetic (Proc. Conf., Battelle Memorial Inst., Seattle, Wash., 1972),  pp. 243--280. Lecture Notes in Math. Vol. 342, Springer, Berlin, 1973. 
%
%\bibitem{Den76}
%R. K. Dennis, {\it In search of new ``homology'' functors having a close relationship
%to K-theory},
%Preprint, Cornell University, Ithaca, NY, 1976.
%
%\bibitem{Eic08}
%B. Eick, and W. Nickel, 
%{Computing the Schur multiplicator and the nonabelian tensor square of a polycyclic group},
%J. Algebra {\bf 320} (2008), 927--944. 
%
%\bibitem{Ell98}
%G. Ellis, {\it The Schur multiplier of a pair of groups},
%Appl. Categ. Struct. {\bf 6} (1998), 355--371.
%
%\bibitem{Ell95}
%G. Ellis, and F. Leonard, {\it Computing Schur multipliers and tensor products of finite groups},
%Proc. Roy. Irish Acad. Sect. A {\bf 95} (1995), no. 2, 137--147.
%
%\bibitem{GAP4}
%The GAP~Group, \emph{GAP -- Groups, Algorithms, and Programming, 
%Version 4.4.12}; 
%2008,
%\verb+(http://www.gap-system.org)+.
%
%\bibitem{Gra79}
%D. R. Grayson, {\it $\K_2$ and the K-theory of automorphisms},
%J. Algebra {\bf 58} (1979), 12--30.
%
%\bibitem{Gru60}
%K. W. Gruenberg, {\it Resolutions by relations},
%J. London Math. Soc. {\bf 35} (1960), 481--494. 
%
\bibitem{Hal40}
P. Hall, {\it The classification of prime-power groups},
J. Reine Angew. Math. {\bf 182} (1940), 130--141.
%
%\bibitem{Hos11}
%A. Hoshi, and M. Kang, {\it Unramified Brauer groups for groups of order $p^5$},
%ArXiv:1109.2966v1, 2011.
%
\bibitem{Hos12}
A. Hoshi, M. Kang, and B. E. Kunyavski\u\i,
{\it Noether problem and unramified Brauer groups},
ArXiv:1202.5812v1, 2012.
%
%\bibitem{Hu09}
%S. Hu, and M. Kang, {\it Noether's problem for some $p$-groups},
%Progress Math., to appear.
%
%\bibitem{Hup67}
%B. Huppert, {\it Endliche Gruppen}, Springer-Verlag, Berlin, 1967.
%
%\bibitem{Iva94}
%S. V. Ivanov, {\it The free Burnside groups of sufficiently large exponents},
%Internat. J. Algebra Comput. {\bf 4} (1994), 1--308.
%
%\bibitem{Jam80}
%R. James, {\it The groups of order $p^6$ ($p$ an odd prime)}, Math. Comp. {\bf 34} (1980), 613--637.
%
%\bibitem{Kar87}
%G. Karpilovsky,
%{\it The Schur multiplier},
%London Mathematical Society Monographs. New Series, 2. The Clarendon Press, Oxford University Press, New York, 1987.
%
\bibitem{Kun08}
B. E. Kunyavski\u\i, {\it The Bogomolov multiplier of finite simple groups},
Cohomological and geometric approaches to rationality problems,  209--217, Progr. Math., 282,
Birkh\"{a}user Boston, Inc., Boston, MA, 2010.
%
%\bibitem{Lyn50}
%R. C. Lyndon, {\it Cohomology theory of groups with a single defining relation},  Ann. of Math. (2)  {\bf 52} (1950), 650--665.
%
%\bibitem{Mag76}
%W. Magnus, A. Karrass, and D. Solitar, {\it Combinatorial group theory},
%Dover Publications, New York, 1976.
%
\bibitem{Mil52}
C. Miller, {\it The second homology of a group},
Proc. Amer. Math. Soc. {\bf 3} (1952), 588--595.
%
%\bibitem{Mil71}
%J. Milnor, {\it Introduction to algebraic K-theory}, Princeton University Press and University of Tokio Press, Princeton, New Jersey, 1971. 
%
%\bibitem{Mor08}
%P. Moravec, {\it On the exponent semigroups of finite p-groups}, J. Group Theory {\bf 11} (2008), no. 4, 511-524.
%
\bibitem{Mor11}
P. Moravec, {\it Unramified groups of finite and infinite groups}, Amer. J. Math., to appear.
%
\bibitem{Mor12}
P. Moravec, {\it Groups of order $p^5$ and their unramified Brauer groups}, submitted.
%
%\bibitem{New68}
%B. B. Newman, {\it Some results on one-relator groups},  Bull. Amer. Math. Soc.  {\bf 74} (1968), 568--571. 
%
\bibitem{Noe16}
E. Noether, {\it Gleichungen mit vorgeschriebener Gruppe}, Math. Ann. {\bf 78} (1916), 221--229. 
%
%\bibitem{Ols91}
%A. Yu. Ol'shanskii, {\it Geometry of defining relations in groups},
%Kluwer Academic Publishers, Dordrecht, Boston, London, 1991.
%
%\bibitem{Pey08}
%E. Peyre, {\it Unramified cohomology of degree 3 and Noether's problem},
%Invent. Math. {\bf 171} (2008), 191--225.
%
%\bibitem{Roc91}
%N. R. Rocco, {\it A construction related to the nonabelian tensor square of a group},
%Bol. Soc. Brasil Mat. (N.S.) {\bf 22} (1991), no. 1, 63--79.
%
\bibitem{Saf91}
I. R. \v{S}afarevi\v{c}, {\it The L\"{u}roth problem},
Proc. Steklov Inst. Math. {\bf 183} (1991), 241--246.
%
\bibitem{Sal84}
D. J. Saltman, {\it Noether's problem over an algebraically closed field},
Invent. Math. {\bf 77} (1984), 71--84.
%
%\bibitem{Sch07}
%I. Schur, {\it Untersuchungen \"{u}ber die Darstellung der endlichen Gruppen durch gebrochene lineare Substitutionen},
%J. Reine Angew. Math. {\bf 132} (1907), 85--137.
%
%\bibitem{Sim94}
%C. C. Sims, {\it Computation with finitely presented groups},
%Encyclopedia of Mathematics and its Applications, vol. 48,
%Cambridge University Press, 1994.
%
%\bibitem{Ste73}
%M. R. Stein, and R. K. Dennis, {\it $\K_2$ of radical ideals and semi-local rings revisited},
%Algebraic $K$-theory, II: ``Classical'' algebraic K-theory and connections with arithmetic (Proc. Conf., Battelle Memorial Inst., Seattle, Wash., 1972),  pp. 281--303. Lecture Notes in Math. Vol. 342, Springer, Berlin, 1973. 
%
\bibitem{Swa83}
R. G. Swan, {\it Noether's problem in Galois theory},
in `Emmy Noether in Bryn Mawr', Springer-Verlag, Berlin, 1983.
%
%\bibitem{Tah72}
%K. Tahara, {\it On the second cohomology groups of semidirect products},
%Math. Z. {\bf 129} (1972), 365--379.
%
\end{thebibliography}
\end{document}